\documentclass[12pt]{amsart}

\newtheorem{theorem}{Theorem}[section]

\usepackage[dvipdf]{graphics}
\usepackage{amssymb}
\usepackage{amsfonts}
\usepackage{latexsym}
\usepackage{cite}
\usepackage{mathrsfs}

\theoremstyle{definition}

\theoremstyle{remark}
\newtheorem{remark}[theorem]{Remark}

\theoremstyle{conjecture}

\theoremstyle{corollary}

\theoremstyle{proposition}
\newtheorem{proposition}[theorem]{Proposition}

\theoremstyle{problem}

\numberwithin{equation}{section}

\large\normalsize

\begin{document}
\begin{center}
{\Large\bf Solving the dimer problem of the vertex-edge graph of a cubic graph}
\\[20pt]
{Shuli\ Li$^{\rm 1}$, Weigen\ Yan$^{\rm 2,*}$ and Danyi Li$^{\rm 2}$} \footnote{ * Corresponding Author.
\newline\hspace*{5mm}{Email address: lishuli198710@163.com} (S. Li), weigenyan@263.net (W. Yan), danyilee1007@163.com (D. Li).
\newline\hspace*{5mm}{The first author was supported by NSFC Grant (11701324) and the Program for Outstanding
Young Scientific Research Talents in Fujian Province University;
The second author was supported by NSFC Grant (12071180; 11571139)}. }
\\[10pt]
\footnotesize { 
$^{\rm 1}$School of Mathematical and Computer Sciences, Quanzhou Normal University,\\ Quanzhou 362000, China
\\[7pt]
$^{\rm 2}$School of Sciences, Jimei University, Xiamen 361021, China}
\\[7pt]

\end{center}
\title{}
\begin{abstract}
Let $G$ be a graph with vertex set $V(G)$ and edge set $E(G)$, and $L(G)$ be the line graph of $G$,
which has vertex set $E(G)$ and two vertices $e$ and $f$ of $L(G)$ is adjacent if $e$ and $f$ is incident in $G$.
The vertex-edge graph $M(G)$ of $G$ has vertex set $V(G)\cup E(G)$ and edge set $E(L(G))\cup \{ue,ve|\ \forall\ e=uv\in E(G)\}$.
In this paper, by a combinatorial technique, we show that if $G$ is a connected cubic graph with an even number of edges,
then the number of dimer coverings of $M(G)$ equals $2^{|V(G)|/2+1}3^{|V(G)|/4}$. As an application, we obtain the exact solution
of the dimer problem of the weighted solicate network obtained from the hexagonal lattice in the context of statistical physics.
\\
{\sl Keywords:}\quad Dimer problem; Vertex-edge graph; Cubic graph; Hexagonal lattice; Solicate network.\\
{\sl MSC(2010):}\quad 05C30; 05C90
\end{abstract}

\maketitle

\section{Introduction}
The graphs considered in this paper are simple, if not specified.
Let $G$ be a connected graph with vertex set $V(G)=\{v_1,v_2,\cdots,v_n\}$ and
edge set $E(G)=\{e_1,e_2,\cdots,e_m\}$. The line graph of $G$, denoted by $L(G)$, is defined as the graph with $V(L(G))=E(G)$ such that
two vertices $e$ and $f$ of $L(G)$ is adjacent if $e$ and $f$ is incident in $G$. The vertex-edge graph of $G$, denoted by $M(G)$, is obtained from $G$ by inserting
a new vertex into each edge of $G$, then joining with edges those pairs of new vertices on adjacent edges of $G$.
Hence it has vertex set $V(G)\cup E(G)$ and edge set $E(L(G))\cup \{ue,ve|\ \forall\ e=uv\in E(G)\}$. For convenience, each vertex $u\in V(G)$ of $M(G)$ is called the $v$-vertex of $M(G)$, and each vertex $f\in E(G)$ of $M(G)$ is called the $e$-vertex of $M(G)$.  The vertex-edge graph of $G$ is also called the middle graph of $G$ in \cite{CDS80}. For the graph $G$ in Figure 1(a), the corresponding line graph $L(G)$ and the vertex-edge graph $M(G)$ are illustrated
in Figure 1(b) and (c), respectively.

\begin{figure}[htbp]
  \centering
\scalebox{1} {\includegraphics{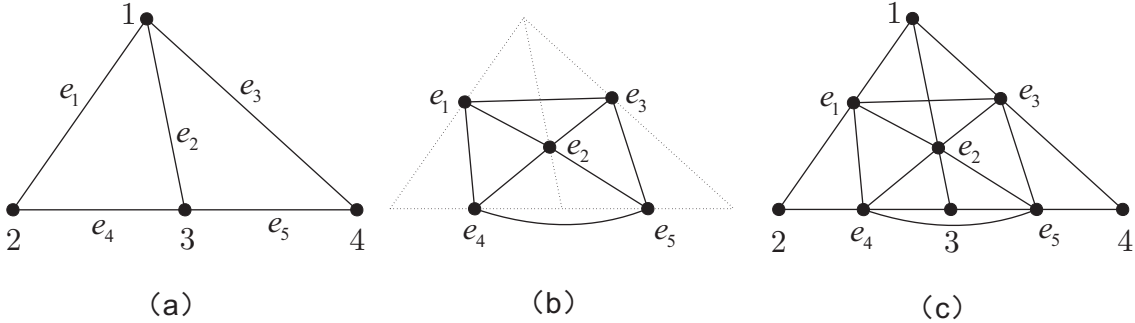}}
  \caption{\ (a)\ A connected graph $G$; (b)\ The line graph $L(G)$ of $G$,
  which is illustrated with solid lines; (c)\ The vertex-edge graph $M(G)$ of $G$.}
\end{figure}

 A set $M\subseteq E(G)$ is a matching of $G$ if every
vertex of $G$ is incident with at most one edge in $M$; it is a perfect matching if every vertex
of $G$ is incident with exactly one edge in $M$. Let $\mathcal{PM}(G)$ be the set of
perfect matchings of $G$ and $p(G)=|\mathcal{PM}(G)|$. Hence $p(G)$ denotes the number of perfect matchings of $G$.
If $G$ is an edge-weighted graph with edge-weight function $\omega: E(G)\rightarrow \mathcal R^+$ ($\mathcal R^+$ is
the set of positive real numbers), the weight of a perfect matching $M$ is the product of weights of edges in $M$.
We also use $p(G)$ to denote the sum of weights of perfect matchings of $G$.

Perfect matchings of a graph are called dimer coverings in statistical physics \cite{Pro99} and
Kekul\'e structures in quantum chemistry \cite{Pau39,Hall73,SHG75}. The dimer models has been widely researched
by mathematicians and physicists \cite{Ciucu97,FL03,Pro99,RST99,Sch98}. It is well known that computing $p(G)$ of a
graph $G$, which is equivalent to the dimer problem in the statistical physics, is an NP-hard problem (see \cite{Jer87,LP86,Val79}).

By the Pfaffian method, the dimer problem on the plane lattices, such as quadratic lattice, 3.12.12 lattice,
Kagom\'e lattice, the Sierpinski gasket, and et al., has been studied extensively by statistical physicists \cite{Wu06,WW07,WW08}. By combinatorial techniques, Ciucu \cite{Ciucu96,Ciucu97,Ciucu98,Ciucu03} obtained the solution of many graphs with some symmetric properties.

For the dimer problem on the line graphs, it is clear that $p(L(G))=0$ if $|E(G)|$ is odd.
Sumner \cite{Sum74} and Vergnas \cite{Ver74} independently showed that every connected claw-free graph with an even number of edges
has one perfect matching. Since every line graph is claw-free, the line graph of a connected graph with an
even size has at least one perfect matching. The enumerative problem of perfect matchings of the line graph
has been independently studied by Ciucu, Liu, Yang \cite{CLY12} and Dong, Yan, Zhang \cite{DYZ13},
they obtained the following result.

\begin{theorem}[\cite{CLY12,DYZ13}]
Let $G=(V(G),E(G))$ be a connected graph, where $|E(G)|$ is even. If the maximum degree $\Delta(G)$ of $G$ satisfies $\Delta(G)\leq 3$, then
\begin{equation}
p(L(G))=2^{|E(G)|-|V(G)|+1}.
\end{equation}
\end{theorem}

Note that the  Kagom\'e, 3.12.12 lattices and the Sierpinski gasket are the line graphs of some graphs, respectively.
It has been determined the number of dimer coverings of the Kagom\'e lattice \cite{WW08}, 3.12.12 lattice \cite{CLY12,DYZ13}, and Sierpinski gasket with
dimension two \cite{DYZ13}. For the vertex-edge graph, the enumerative problem of spanning trees \cite{Yan17},
the normalised Laplacian characteristic polynomial and degree-Kirchhoff index \cite{HL15} has been studied.
A natural problem is how to solve the dimer problem of the vertex-edge graph of a graph.

{\bf The one of the main purposes} in this paper is to obtain the solution of the dimer problem of the
vertex-edge graph of the graphs whose maximum degree is no more than three.

The following result is obvious.

\begin{proposition}
If a graph $G$ has a pendent vertex $v$, then $p(M(G))=p(M(G-v))$.
\end{proposition}

If a graph $H$ has one vertex $u$ of degree 2, and $x$ and $y$ are the two adjacent vertices of $u$.
Let $H'$ be the graph obtained from $H$ by deleting vertex $u$ and then joining an edge $xy$.
Then we have the following result.

\begin{proposition}
Suppose $H$ and $H'$ are the graphs defined above. Then
$$p(M(H))=p(M(H')).$$
\end{proposition}
\begin{proof}
The result is a direct corollary of Lemma 1.3 in \cite{Ciucu97}.
\end{proof}

Suppose that $H$ is a connected graph whose maximum degree $\Delta(H)\leq 3$ and $|E(H)|$ is even. We know that the number of vertices of $H'$ of degree 2 is one less than that in $H$. By repeating the two operations in Propositions 1.2 and 1.3, we can obtained a graph $H^*$ which is a cubic graph or a cycle with two vertices. Particularly, $p(M(H))=p(M(H^*))$. Hence the dimer problem on the vertex-edge graph $M(G)$ of a graph $G$ with the maximum degree no more than three is equivalent to the dimer problem on cubic graphs. So in the following, we only need to consider the dimer problem of the vertex-edge graph of cubic graphs.

 We will prove the following results in the next section.

\begin{theorem}
Let $G$ be a connected cubic graph with $n$ vertices and $m$ edges, where $m$ is even. Then
\begin{equation}
p(M(G))=2^{m-n+1}3^{\frac{2n-m}{2}}=2^{\frac{n}{2}+1}3^{\frac{n}{4}}.
\end{equation}
\end{theorem}

\begin{theorem}
Let $G$ be a connected cubic graph with $n$ vertices and $m$ edges and $e$ be a non-cut edge of $G$, where $m$ is odd. Then
\begin{equation}
p(M(G-e))=p(M(G)-e)=2^{m-n}3^{\frac{2n-m-1}{2}}=2^{\frac{n}{2}}3^{\frac{n-2}{4}}.
\end{equation}
\end{theorem}

\begin{theorem}
Let $G$ be a connected cubic graph with $n$ vertices and $m$ edges and $e$ be a cut edge of $G$, where $m$ is odd. Suppose that $G-e=G_1\cup G_2$.

(a). \ If both $G_1$ and $G_2$ have an even number of edges, then
$$p(M(G-e))=p(M(G)-e)=0.$$

(b). \  If both $G_1$ and $G_2$ have an odd number of edges, then
\begin{equation}
p(M(G-e))=p(M(G)-e)=2^{m-n+1}3^{\frac{2n-m-1}{2}}=2^{\frac{n}{2}+1}3^{\frac{n-2}{4}}.
\end{equation}
\end{theorem}

{\bf As an application}, in Section 3 we obtain the exact solution
of the dimer problem on the weighted solicate network obtained from the hexagonal lattice in the context of statistical physics. Finally, in the last section, we give some discussion.

\section{Proofs of main results}

First, we give a combinatorial proof of Theorem 1.4 as follows.

{\bf Proof of Theorem 1.4}.\ \  Note that $G$ is a cubic graph with $n$ vertices and $m$ edges. So $m=\frac{3n}{2}$. Since $m$ is even, $4|n$. For each vertex $v_{k}\in V(G), k=1,2,\ldots,n$,
if we denote the three edges incident with $v_k$ in $G$ by $e_{k_1},e_{k_2}$ and $e_{k_3}$,
then the four vertices $v_{k},e_{k_1},e_{k_2},e_{k_3}\in V(M(G))$ induce a complete subgraph $K_4$ in $M(G)$ ($K_4$ is a complete graph with four vertices), denoted by $G_k$. Obviously, $G_k$ has three $e$-vertices $e_{k_1},e_{k_2},e_{k_3}$ and one $v$-vertex $v_k$. Thus the vertex-edge graph $M(G)$ can be decomposed into $n$ edge-disjoint $K_4$'s $G_1, G_2,\ldots,G_n$, i.e., $E(M(G))=E(G_1)\cup E(G_2)\cup \cdots\cup E(G_n)$.
Note that $G$ is cubic and $m$ is even. Hence, by Theorem 1.1,
\begin{equation}
p(L(G))=2^{m-n+1}=2^{\frac{n}{2}+1}.
\end{equation}

\begin{figure}[htbp]
  \centering
\scalebox{0.9}{\includegraphics{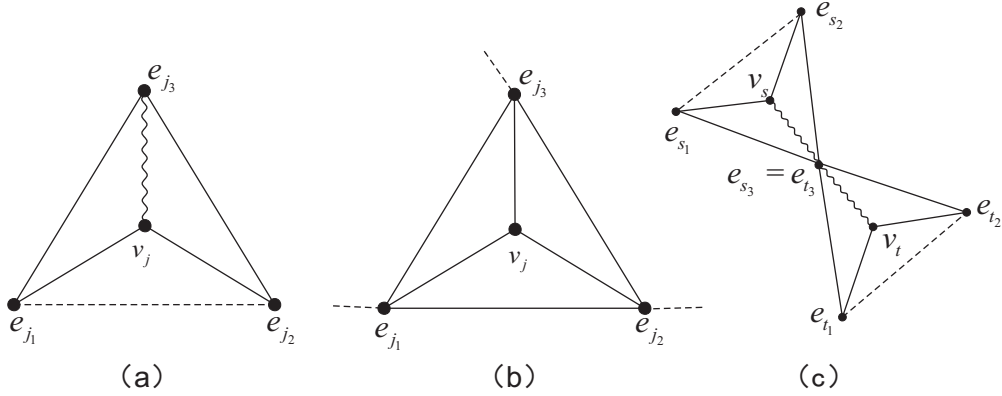}}
  \caption{\ (a)\ The $K_4$ subgraph $G_j$ of $M(G)$ with dotted edge $e_{j_1}e_{j_2}\in M_l$, where $M_l$ is a dimer covering of $L(G)$;
  (b)\ The $K_4$ subgraph $G_j$ of $M(G)$ such that $E(G_j)\cap M_l=\emptyset$, i.e., the three dotted edges are in $M_l$.\ (c)\ Two graphs $G_s$ and $G_t$ with one common vertex.}
\end{figure}

Let $\mathcal{PM}(L(G))$ be the set of dimer coverings of $L(G)$. Since $p(L(G))=|\mathcal{PM}(L(G))|=2^{m-n+1}=2^{\frac{n}{2}+1}$, we may set $\mathcal {PM}(L(G))=\{M_1,M_2,\ldots, M_{2^{\frac{n}{2}+1}}\}$.

For any fixed dimer covering $M_{l}\in \mathcal{PM}(L(G))=\{M_i|1\leq i\leq 2^{\frac{n}{2}+1}\}$, $|M_{l}|=\frac{m}{2}=\frac{3n}{4}$. Obviously, $M_l$ is a matching of $M(G)$. Particularly, we have the following claims.

{\bf Claim 1}.\ \ There exist exactly $\frac{m}{2}=\frac{3n}{4}$ edge-disjoint $K_4$'s $G_{l_1},G_{l_2}, \ldots, G_{l_{\frac{3n}{4}}}$ in $\{G_i|1\leq i\leq n\}$ such that for any $j=l_1,l_2,\ldots,l_{\frac{3n}{4}}$, $G_{j}$ contains exactly one edge in $M_l$ (see the dotted edge $e_{j_1}e_{j_2}$ in Figure 2(a)).

For each $j\in \{l_1,\l_2,\ldots,l_{\frac{3n}{4}}\}$, $G_j$ has three $e$-vertices $e_{j_1},e_{j_2}, e_{j_3}$ and one $v$-vertex $v_j$. Without loss of generality, we assume that edge $e_{j_1}e_{j_2}$ of $G_j$ is in $ M_l$ (see Figure 2(a)). That is, $M_l=\{e_{j_1}e_{j_2}|j=l_1,l_2,\ldots,l_{\frac{3n}{4}}\}$.

{\bf Claim 2}.\ \ For each $G_j, j\in \{1,2,\ldots,n\}\setminus \{l_1,\l_2,\ldots,l_{\frac{3n}{4}}\}$,  the three $e$-vertices $e_{j_1},e_{j_2},e_{j_3}$ of $G_j$ are matched by three different edges in $M_l$ non of which is in $E(G_j)$ (see Figure 2(b)).

{\bf Claim 3}.\ \ The $n-\frac{m}{2}=n-\frac{3n}{4}=\frac{n}{4}$ $G_j$'s ($j\in \{1,2,\ldots,n\}\setminus \{l_1,\l_2,\ldots,l_{\frac{3n}{4}}\}$) are vertex-disjoint complete subgraphs of $M(G)$ with four vertices.

Set $I=\{1,2,\ldots,n\}\setminus \{l_1,\l_2,\ldots,l_{\frac{3n}{4}}\}$ and $G_l'=\cup_{s\in I}G_s$.
Let $\mathcal M_l'$ be the set of dimer coverings of $G_l'$. By Claim 3,
\begin{equation}
|\mathcal M_l'|=3^{\frac{n}{4}},
\end{equation}
since each $G_s$ ($s\in I$) has three dimer coverings.

Now, from each dimer covering $M_l$ of $L(G)$, $ l\in \{1,2,\ldots,2^{m-n+1}=2^{\frac{n}{2}+1}\}$, we construct a subset $\mathcal M_l^* \subseteq \mathcal{PM}(M(G))$ as follows.

By Claim 1, for any $j\in \{l_1,l_2,\ldots,l_{\frac{3n}{4}}\}$, each $G_{j}$ has exactly one edge $e_{j_1}e_{j_2}$ in $M_l$. Delete this edge and its incident two vertices from $G_{j}$. Then the resulting subgraph is an edge $v_je_{j_3}$ which has one $v$-vertex $v_j$ and one $e$-vertex $e_{j_3}$ of $M(G)$ (see the wave edge in Figure 2(a)).

Set $$M'_l=\{v_je_{j_3}|j=l_1,l_2,\ldots,l_{\frac{3n}{4}}\}$$
and
\begin{equation}
\mathcal M_{l}^{*}=\{M'_{l}\cup M|M\in \mathcal M_l'\}.
\end{equation}

{\bf Claim 4}.\ \ $M_l'=\{v_je_{j_3}|j=l_1,l_2,\ldots,l_{\frac{3n}{4}}\}$  is a matching of $M(G)$.

Otherwise, there exist two edges $v_se_{s_3}$ and $v_te_{t_3}$ in $M_l'$ such that $e_{s_3}=e_{t_3}$, where $s, t\in \{l_1,l_2,\ldots,l_{\frac{3n}{4}}\}$ and $s\neq t$. This implies that $e_{s_3}=e_{t_3}$ which is an edge of $G$ with form of $v_sv_t$ and hence $G_s$ and $G_t$ have a common vertex $e_{s_3}=e_{t_3}$ (see Figure 2(c)). By the definition of $M_l$, $e_{s_1}e_{s_2}$ and $e_{t_1}e_{t_2}$ are two edge in $M_l$. Then the vertex $e_{s_3}=e_{t_3}$ of the line graph $L(G)$ can not be matched by the edges in $M_l$. This is a contradiction since $M_l$ is a dimer covering of $L(G)$. Hence the claim holds.

Obviously, by Eqs. (2.2) and (2.3),
\begin{equation}
|\mathcal M_{l}^{*}|=3^{\frac{n}{4}}.
\end{equation}
Moreover, we have the following claims.

{\bf Claim 5}.\ \ $\mathcal M_{l}^{*}\subseteq \mathcal {PM}(M(G))$. That is, for any $M\in \mathcal M_l'$, $M'_{l}\cup M$ is a dimer covering of $M(G)$.

By Claim 3, $M$ is a dimer coverings of $G_l'=\cup_{s\in I}G_s$, which contains exactly $2(n-\frac{m}{2})=2(n-\frac{3n}{4})=\frac{n}{2}$ vertex-disjoint edges, which is a matching of $M(G)$. Since $|M_l'\cup M|=|M_l'|+|M|=\frac{3n}{4}+\frac{n}{2}=\frac{5n}{4}=\frac{n+m}{2}$, it suffices to prove that $M_l'\cup M$ is a matching of $M(G)$ (hence it is a perfect matching or a dimer covering of $M(G)$). By the definition of $\mathcal M_l^*$, $M$ is a matching of $M(G)$. By Claim 4, $M_l'$ is a matching of $M(G)$. Hence we only need to prove that any edge in $M$ and any edge in $M_l'$  have no common vertex. For each $s\in I=\{1,2,\ldots,n\}\setminus \{l_1,l_2,\ldots,l_{\frac{3n}{4}}\}$, $G_s$ has three dimer coverings (perfect matchings): $\{v_se_{s_1},e_{s_2}e_{s_3}\}, \{v_se_{s_2},e_{s_1}e_{s_3}\}, \{v_se_{s_3},e_{s_1}e_{s_1}\}$. By Claim 2, the three $e$-vertices $e_{s_1}, e_{s_2}, e_{s_3}$ of $G_s$ are matched by three different edges in $M_l$ (see Figure 2(b)). By the definition of $M_l'$, any edge in $M$ and any edge in $M_l'$ have no common vertex. Hence Claim 5 holds.

{\bf Claim 6}.\ \ For any two different dimer coverings $M_p$ and $M_q$ of $L(G)$, $\mathcal M_p^*\cap \mathcal M_q^*=\emptyset$.

By the definition of $\mathcal{M}_l^*$, Claim 6 is obvious.

{\bf Claim 7}.
\begin{equation}
\cup_{l=1}^{2^{\frac{n}{2}}+1}\mathcal M_l^*=\mathcal {PM}(M(G)).
\end{equation}

By Claim 5, $\cup_{l=1}^{2^{\frac{n}{2}}+1}\mathcal M_l^*\subseteq \mathcal {PM}(M(G))$. It suffices to prove that $\cup_{l=1}^{2^{\frac{n}{2}}+1}\mathcal M_l^*\supseteq\mathcal {PM}(M(G))$. Given a dimer covering $M^*$ of $M(G)$, $|M^*|=\frac{5n}{4}$. Since $M(G)$ contains $n$ $v$-vertices $v_1,v_2,\ldots,v_n$ and $\frac{3n}{2}$ $e$-vertices, $M^*$ contains $\frac{1}{2}\left(\frac{3n}{2}-n\right)=\frac{n}{4}$ edges each of which is incident two $e$-vertices  and $n$ edges each of which is incident one $v$-vertex and one $e$-vertex. Without loss of generality, denote the $n$ edges in $M^*$ each of which is incident one $v$-vertex and one $e$-vertex by $v_1e_{1_3},v_2e_{2_3},\ldots,v_ne_{n_3}$, respectively. Assume that there $x$ graphs $G_{a_1},G_{a_2},\ldots,G_{a_x}$ in $\{G_1,G_2,\ldots,G_n\}$ each of which has two edges in $M^*$, and $y$ graphs $G_{b_1},G_{b_2},\ldots,G_{b_y}$ in $\{G_1,G_2,\ldots,G_n\}$ each of which has one edges in $M^*$. Then $x+y=n, 4x+2y=n+m=\frac{5n}{2}$. Hence $x=\frac{n}{4}, y=\frac{3n}{4}$. Note that, for any $t\in \{b_1,b_2,\ldots,b_{\frac{3n}{4}}\}$, $G_t$ has the vertex set $\{v_t,e_{t_1},e_{t_2},e_{t_3}\}$ and the edge set $\{v_te_{t_1},v_te_{t_2},v_te_{t_3},e_{t_1}e_{t_2},e_{t_2}e_{t_3},e_{t_3}e_{t_1}\}$.  The unique edge of $G_{t}$ which belongs to $M^*$ has one $v$-vertex $v_{t}$ and one $e$-vertex, say $e_{t_3}$. This implies that $\{v_te_{t_3}|t=b_1,b_2,\ldots,b_{\frac{3n}{4}}\}\subseteq M^*$.

{\bf Claim 7.1}.\ \ $\{e_{t_1}e_{t_2}|t=b_1,b_2,\ldots,b_{\frac{3n}{4}}\}$ is a dimer covering (perfect matching) of $L(G)$ (the line graph of $G$).

Otherwise, there exist $p,q\in \{b_1,b_2,\ldots,b_{\frac{3n}{4}}\}$ ($p\neq q$) such that two edges $e_{p_1}e_{p_2}$ and $e_{q_1}e_{q_2}$ have a common vertex, say $e_{p_1}=e_{q_1}$. That is, $G_p$ and $G_q$ have a common vertex $e_{p_1}=e_{q_1}$. Note that there exists a unique edge $v_pe_{p_3}$ of $G_p$ which belongs to $M^*$ (resp. a unique edge $v_qe_{q_3}$ of $G_q$ which belongs to $M^*$). This implies that the vertex $e_{p_1}=e_{q_1}$ of $M(G)$ can not be matched by $M^*$, a contradiction (since $M^*$ is a dimer covering of $M(G)$). Hence Claim 7.1 holds.

By Claim 7.1, $M_k=\{e_{t_1}e_{t_2}|t=b_1,b_2,\ldots,b_{\frac{3n}{4}}\}\subseteq \mathcal{PM}(L(G))$. Let $\mathcal M_k^*$ be the set of dimer coverings of $M(G)$ defined in
Eq. (2.3). It is not difficult to see that $M^*\in \mathcal M_k^*$, which implies that $\cup_{l=1}^{2^{\frac{n}{2}}+1}\mathcal M_l^*\supseteq\mathcal {PM}(M(G))$. Claim 7 thus follows.

Using Claims 6 and 7 and Eq. (2.4), then
$$p(M(G))=|\mathcal {PM}(M(G))|=|\cup_{l=1}^{2^{\frac{n}{2}+1}}\mathcal M_l^*|=\sum_{l=1}^{2^{\frac{n}{2}+1}}|\mathcal M_l^*|=2^{\frac{n}{2}+1}3^{\frac{n}{4}}. \quad\Box $$

Now, we can prove Theorem 1.5.

{\bf Proof of Theorem 1.5}.\ \ Note that $G$ is a cubic graph. Then $3n=2m$.
 Since $m$ is odd and $n$ is even, $|V(M(G))|=n+m$ is odd. Note that $G-e$ is a connected graph with two vertices of degree 2. Using two operations in Proposition 1.3 to $G-e$, we obtain a
connected cubic graph $G'=(V(G'),E(G'))$ with $n-2$ vertices and $m-3$ edges such that $p(M(G-e))=p(M(G'))$.
By Theorem 1.4,
\begin{align*}
p(M(G-e))
&=p(M(G'))\\
&=2^{m-3-(n-2)+1}3^{\frac{2(n-2)-(m-3)}{2}}\\
&=2^{\frac{n}{2}}3^{\frac{n-2}{4}}.
\end{align*}

Hence the theorem follows.

{\bf Proof of Theorem 1.6}. Note that $e=uv$ is a cut edge of $G$, $G-e=G_1\cup G_2$, and $V(G_1)\cap V(G_2)=\emptyset$. Hence
\begin{equation}
p(M(G-e))=p(M(G_1))p(M(G_2)).
\end{equation}
Set $|V(G_1)|=n_1, |E(G_1)|=m_1$, $|V(G_2)|=n_2$, $|E(G_2)|=m_2$.
Then $n_1+n_2=n, m_1+m_2+1=m$. For the graph $G_1$ (resp. $G_2$), there is exactly a vertex $u$ (resp. $v$) of degree 2. Then
\begin{equation}
3(n_1-1)+2=2m_1, 3(n_2-1)+2=2m_2.
\end{equation}
By Eq. (2.7), both of $n_1$ and $n_2$ are odd. Note that $M(G_1)$ (resp. $M(G_2)$) has $n_1+m_1$ (resp. $n_2+m_2$) vertices,
and $m=m_1+m_2+1$ is odd. So $m_1$ and $m_2$ have the same parity.
If $m_1$ (resp. $m_2$) is even, then $n_1+m_1$ (resp. $n_2+m_2$) is odd, so $p(M(G_1))=0$ (resp. $p(M(G_2))=0$), then by Eq. (2.6),
$p(M(G-e))=0$. So in the following we suppose that $m_1$ (resp. $m_2$) is odd.

Note that both $G_1$ and $G_2$ have one vertex of degree two. Using the operation in Proposition 1.3, we can obtained two cubic graph $G'_1$ and $G'_2$ with $|V(G'_1)|=n_{1}-1, |E(G'_1)|=m_{1}-1$ and $|V(G'_2)|=n_{2}-1, |E(G'_2)|=m_{2}-1$. Particularly,
\begin{equation}
p(M(G_1))=p(M(G'_1)), p(M(G_2))=p(M(G'_2)).
\end{equation}
Note that $|E(G'_1)|=m_{1}-1$ (resp. $|E(G'_2)|=m_{2}-1$) is even. Hence, by Theorem 1.4,
\begin{equation}
p(M(G'_1))
=2^{(m_1-1)-(n_1-1)+1}3^{\frac{2(n_1-1)-(m_1-1)}{2}}
=2^{m_1-n_1+1}3^{\frac{2n_1-m_1-1}{2}}.
\end{equation}
Similarly,
\begin{equation}
p(M(G'_2))=2^{m_2-n_2+1}3^{\frac{2n_2-m_2-1}{2}}.
\end{equation}
By Eqs. (2.6), (2.8)-(2.10), we have
\begin{align*}
p(M(G-e))
&=p(M(G'_1))p(M(G'_2))\\
&=2^{m_1-n_1+1}3^{\frac{2n_1-m_1-1}{2}}2^{m_2-n_2+1}3^{\frac{2n_2-m_2-1}{2}}\\
&=2^{m-n+1}3^{\frac{2n-m-1}{2}}\\
&=2^{\frac{n+2}{2}}3^{\frac{n-2}{4}}.
\end{align*}
The theorem thus follows.

\section{An application}

In this section, as an application of Theorem 1.4, we enumerate dimer coverings of the weighted silicate
networks with toroidal boundary condition which can be obtained from the hexagonal lattices in the context of statistical physics as follows.
\begin{figure}[htbp]
  \centering
\scalebox{1} {\includegraphics{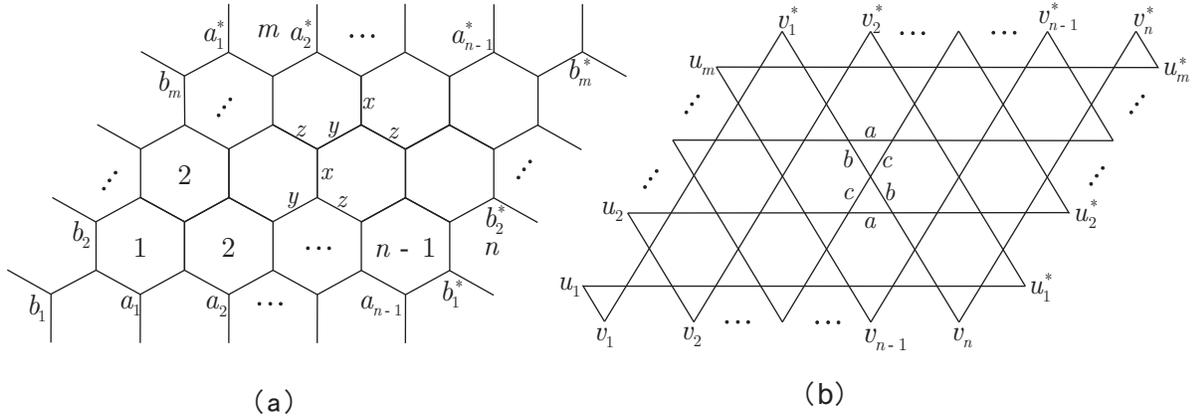}}
  \caption{\ (a)\ The hexagonal lattice $H^{T}(n,m)$ with toroidal boundary condition;
  (b)\ A graph $G(n,m)$. }
\end{figure}

We first introduce the hexagonal lattices which have been extensively studied by statistical physicists \cite{Els89,Fis66,PW88}.
The hexagonal lattice $H^{T}(n,m)$ (where $mn$ is even) with toroidal boundary condition is shown in Figure 3(a),
where $(b_1, b_{1}^*),(b_2,b_{2}^*),\ldots, (b_{m},b_{m}^*)$ and $(b_1,a_{1}^*)$,$(a_1,a_{2}^*),\ldots,$ $(a_{n-2},a_{n-1}^*)$,$(a_{n-1}, b_m^*)$
are edges in $H^{T}(n,m)$.  The edges in $H^{T}(n,m)$ are drawn in three different directions,
we weight each edge on the vertical direction with $x$, each edge on the slash direction with $y$,
each edge on the backslash direction with $z$, as shown in Figure 3(a). Note that $H^{T}(n,m)$ is a cubic graph
with $2mn$ vertices and $3mn$ edges and hence the number of edges in each of the three directions is $mn$.

Let $G(n,m)$ be the plane lattice graph illustrated in Figure 3(b), each of whose vertices has degree two or four.
For $G(n,m)$, if we identify each pair of vertices $u_i$ and $u_{i}^*$ for all $i=1,2,\ldots,m$, and $v_j$ and $v_{j}^*$ for all
$j=1,2,\ldots,n$, the resulting graph, denoted by $K^{T}(n,m)$, is called the Kagom\'e lattice with toroidal boundary
condition by statistical physicists (see \cite{Els89,PW88,WW07,Wu06,WW08}. By the definition of $K^{T}(n,m)$, we know that $K^{T}(n,m)$
has $3mn$ vertices. It is not difficult to see that $K^{T}(n,m)$ is the line graph of $H^{T}(n,m)$. Dong, Yan and Zhang \cite{DYZ13}
proved that if $mn$ is even, then
\begin{equation}
p(K^{T}(n,m))=2^{mn+1}.
\end{equation}
The edges in $K^{T}(n,m)$ are drawn in three different directions as shown in Figure 3(b).
If we weight each edge on the horizontal direction by $a$, each edge on the slash direction by $c$, and each edge on the
backslash direction by $b$, as shown in Figure 3(b), the resulting graph, denoted by $K_{\omega}^{T}(n,m)$.
Wu and Wang \cite{WW08} considered the problem of weighted enumeration of dimer coverings of $K_{\omega}^{T}(n,m)$.
They proved that if $mn$ is even, then the weighted enumeration of dimer coverings of $K_{\omega}^{T}(n,m)$ can be expressed by
\begin{equation}
p(K_{\omega}^{T}(n,m))=2^{mn+1}(abc)^{\frac{mn}{2}}.
\end{equation}

\begin{figure}[htbp]
  \centering
\scalebox{1} {\includegraphics{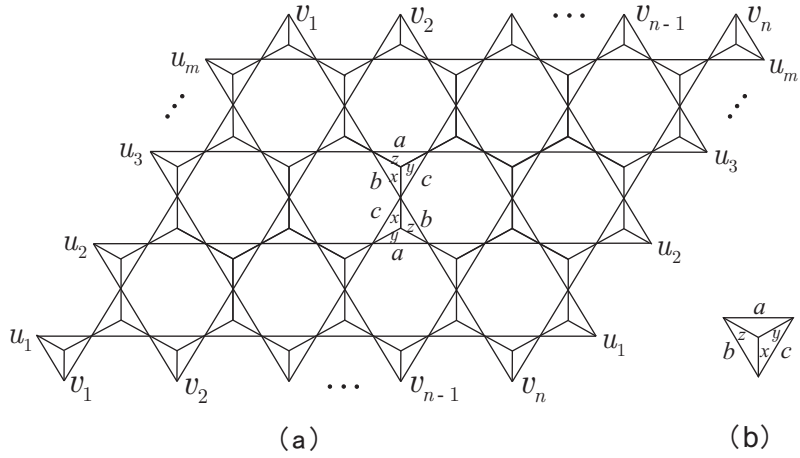}}
  \caption{\ (a).\ The silicate networks $S^{T}(n,m)$. \ (b).\ A weighted $K_4^{\omega}$ with six edges of weights $a,b,c,x,y,z$.}
\end{figure}

Let $S^T(n,m)$ be the lattice on the torus, which is illustrated in Figure 4(a). $S^T(n,m)$ is a kind of the silicate networks with toroidal boundary condition
by statistical physicists and chemists. Solicate networks has been widely studied by physicists, chemists and mathematicians
\cite{MR09,MRWK14,NSZ02,RWGS12,SPX20,ZPIBA12}. By the definition of $S^{T}(n,m)$, we know that $S^{T}(n,m)$ has $5mn$ vertices.
It is not difficult to see that $S^{T}(n,m)$ is the vertex-edge graph (middle graph) of $H^{T}(n,m)$. The edges in $S^{T}(n,m)$ are drawn
in six different directions as shown in Figure 4(a) or (b). If we weight edges of $S^{T}(n,m)$ by six kind of weights $a,b,c,x,y,z$ such that edges with the same direction have the same weight as shown in Figure 4(a) and (b), then
we obtain a weighted silicate networks, denoted by $S_{\omega}^{T}(n,m)$.

\begin{theorem}
Let $S_{\omega}^{T}(n,m)$ be the weighted silicate networks with toroidal boundary condition, where $mn$ is even. Then the sum of weights of dimer coverings of $S_{\omega}^{T}(n,m)$ can be expressed by
$$p(S_{\omega}^{T}(n,m))=2^{mn+1}(xyz)^{\frac{mn}{2}}(ax+by+cz)^{\frac{mn}{2}}.$$
\end{theorem}
\begin{proof}
Note that $K^{T}(n,m)$ is the line graph of $H^{T}(n,m)$. By Eq. (3.2), the sum of weights of dimer coverings of the weighted Kagom\'e lattice $K_{\omega}^{T}(n,m)$ is
$$p(K_{\omega}^{T}(n,m))=p\big(L(H^{T}(n,m))\big)=2^{mn+1}(abc)^{\frac{mn}{2}}.$$
This implies that $K_{\omega}^{T}(n,m)$ has $2^{mn+1}$ dimer coverings and each dimer covering of $K_{\omega}^{T}(n,m)$ has $\frac{mn}{2}$ edges of weight $a$, $\frac{mn}{2}$ edges of weight $b$, and $\frac{mn}{2}$ edges of weight $c$, respectively.

On the other hand, we know that the silicate network $S^T(n,m)$ is the vertex-edge graph of $H^T(n,m)$ and the weighted silicate network $S_{\omega}^{T}(n,m)$ has $2mn$ edge-disjoint $K_4^{\omega}$ illustrated in Figure 4(b). The proof of Theorem 1.4 implies that
$$\mathcal{PM}(S^T(n,m))=\cup_{l=1}^{2^{mn+1}}\mathcal M_l^*,$$
where $\mathcal{PM}(K^T(n,m))=\mathcal{PM}(L(H^T(n,m)))=\{M_l|1\leq l\leq 2^{mn+1}\}$. Since each $M_l\in \mathcal{PM}(K^T(n,m))$, $M_l$ has $\frac{mn}{2}$ edges of weight $a$, $\frac{mn}{2}$ edges of weight $b$, and $\frac{mn}{2}$ edges of weight $c$, respectively.
By Eq. (2.3) in the proof of Theorem 1.4,
$$\mathcal M_{l}^{*}=\{M'_{l}\cup M|M\in \mathcal M_l'\}, $$
where $\mathcal M_l'$ is the set of dimer coverings of the union of $2mn-\frac{3mn}{2}=\frac{mn}{2}$ vertex-disjoint $K_4^{\omega}$. By the definition of $M_l'$ in the proof of Theorem 1.4, $M_l'$ is obtained from $M_l$ by replacing each edge of $M_l$ of weight $a$ with an edge of weight $x$, replacing each edge of weight $b$ with an edge of weight $y$, and replacing each edge of weight $c$ with an edge of weight $z$ (See Figure 4(a) and (b)). Note that the product of weights of edges in $M_l$ equals $(abc)^{\frac{mn}{2}}$. Hence the product of weights of edges in $M_l'$ equals $(xyz)^{\frac{mn}{2}}$. Since the sum of weights of dimer coverings of each $K_4^{\omega}$ equals $ax+by+cz$ (that is $p(K_4^{\omega})=ax+by+cz$), the sum of weights of dimer coverings in $\mathcal M_l^*$ can be expressed by
$$\omega(\mathcal M_l^*)=(xyz)^{\frac{mn}{2}}(ax+by+cz)^{\frac{mn}{2}}.$$
Hence the sum of weights of dimer coverings of $S_{\omega}^T(n,m)$ is
$$p(S_{\omega}^T(n,m))=\sum_{l=1}^{2^{mn+1}}\omega(\mathcal M_l^*)=2^{mn+1}(xyz)^{\frac{mn}{2}}(ax+by+cz)^{\frac{mn}{2}}. $$
So we have finished the proof of Theorem 3.1.
\end{proof}

\begin{remark}
If each edge of $S_{\omega}^{T}(n,m)$ has weight 1, then
$$p(S^{T}(n,m))=2^{mn+1}3^{mn/2},$$
which coincides with the conclusion as Theorem 1.4.
\end{remark}

\section{Discussion}
In this paper, we obtain the solution of the dimer problem on the vertex-edge graphs of cubic graphs. As an application, we give a simple formula of the sum of weights of dimer covering of the weighted silicate network obtained from the hexagonal lattice with toroidal boundary condition in the context of statistical physics. Note that the cubic graphs are a type of important graphs not only in combinatorics but also in physics. For example, $3.12.12, 4.8.8, 4.6.12$, hexagonal lattices, and so on, which are cubic graphs, are very important in statistical physics on which the dimer problem has been extensively studied by statistical physicists (see for example \cite{Ciucu03,Wu06}). The main results in this paper can be applied to solve the dimer problem of the weighted vertex-edge graphs of these cubic graphs.

\end{document}